\documentclass{article}
\usepackage{amssymb,amsthm,graphics}
\newtheorem{theorem}{Theorem}[section]
\newtheorem{lemma}[theorem]{Lemma}

\newcommand{\bq}{\[}
\newcommand{\eq}{\]}
\newcommand{\bqn}{\begin{equation}}
\newcommand{\eqn}{\end{equation}}

\newcommand{\TM}{{${}^{\scriptscriptstyle\rm TM}$}}
\newcommand{\R}{\mathbb{R}}
\newcommand{\Z}{\mathbb{Z}}
\newcommand{\st}{{\bigm|}}

\newcommand{\Det}{\mathrm{Det}}
\newcommand{\Per}{\mathrm{Per}}
\newcommand{\Pf}{\mathrm{Pf}}
\newcommand{\Hf}{\mathrm{Hf}}
\newcommand{\N}{{$\nabla\:$}}
\newcommand{\D}{{$\Delta\:$}}
\newcommand{\T}[1]{\tilde{#1}}
\renewcommand{\O}[1]{\overline{#1}}

\date{ }

\begin{document}

\title{Symmetries of plane partitions and the permanent-determinant method}
\author{Greg Kuperberg}

\maketitle

\begin{abstract}
In the paper \cite{Stanley:symmetries}, Stanley gives formulas for the
number of plane partitions in each of ten symmetry classes.  This paper
together with the two papers \cite{Andrews:tsscpp} and
\cite{Stembridge:enumeration} completes the project of proving all ten
formulas.

We enumerate cyclically symmetric, self-complementary plane
partitions.  We first convert plane partitions to tilings of a hexagon
in the plane by rhombuses, or equivalently to matchings in a certain
planar graph.  We can then use the permanent-determinant method or a
variant, the Hafnian-Pfaffian method, to obtain the answer as the
determinant or Pfaffian of a matrix in each of the ten cases.  We
row-reduce the resulting matrix in the case under consideration to
prove the formula.  A similar row-reduction process can be carried out
in many of the other cases, and we analyze three other symmetry classes
of plane partitions for comparison.
\end{abstract}

\section{Introduction}

In the paper \cite{Stanley:symmetries}, Richard Stanley describes ten
symmetry classes of plane partitions and gives formulas for their
enumeration.  We give a brief summary of the conjectures and results
mentioned in the paper:

A plane partition is a finite union $\pi$ of unit cubes in the positive
octant of $\R^3$ such that for each cube $C$ in $\pi$, there is either
another cube or a wall of the octant below $C$, and the same is true
behind and to the left of $C$.  (Sometimes this arrangement is called
the Ferrer's diagram of a plane partition.) Alternatively, a plane
partition is a finite subset of $\Z_+ \times \Z_+ \times \Z_+$ which is
an order ideal with respect to the partial ordering
$(x+1,y,z),(x,y+1,z),(x,y,z+1) \ge (x,y,z)$ for $(x,y,z) \in \Z_+^3$.
An example of an ordinary but useful plane partition is a box $B(a,b,c)
= [0, a] \times [0,b] \times [0,c]$.  We  define three actions on the
set of pairs $(\pi,B(a,b,c))$, where $\pi \subseteq B(a,b,c)$:

\begin{description}
\item[1.] We define $\tau(\pi)$, the transpose of $\pi$, to be
$\{(x,y,z) \st (y,x,z) \in \pi\}$, and $\tau(\pi,B(a,b,c))$ to be
$(\tau(\pi),B(b,a,c))$.

\item[2.] We define $\rho(\pi)$, the rotation of $\pi$, to be
$\{(x,y,z) \st (y,z,x) \in \pi\}$, and $\rho(\pi,B(a,b,c))$ to be
$(\rho(\pi),B(c,a,b))$.

\item[3.] We define $\kappa(\pi,B(a,b,c)) = (\pi', B(a,b,c))$, where
$\pi'$ is the complement of $\pi$ in $B(a,b,c)$ with
all three coordinates reversed.  Whenever a box $B(a,b,c)$ has been
chosen, we may write $\kappa(\pi)$ for $\pi'$ by abuse of notation.
\end{description}

Let $T$ be the group generated by $\tau$, $\rho$, and $\kappa$.  For
each subgroup $G$ of $T$, we define $N_G(a,b,c)$ to be the number of
pairs $(\pi,B(a,b,c))$ fixed by $G$ (in general, $N_G(a,b,c)$ is zero
unless $B(a,b,c)$ is fixed by $G$).  There are ten inequivalent choices
$G_1,\ldots,G_{10}$ for $G$; we set $N_i(a,b,c) = N_{G_i}(a,b,c)$.
Table~\ref{formulas} gives the definition of each of the ten groups,
together with the corresponding formula for $N_i$.  In the table,
$H(n)$ is the hyperfactorial function and $H_k(n)$ is a staggered
hyperfactorial, defined by $H_k(n) = (n-k)!(n-2k)!(n-3k)!\ldots$,
$H(n) = H_1(n)$, and $F_k(n) = n(n-k)(n-2k)\ldots$ is a staggered
factorial.  The symbol $\stackrel{!}{=}$ indicates that the equality
was conjectural before this paper.

\begin{table} \renewcommand{\arraystretch}{1.5}
$$ \begin{array}{|c|rcl|}
\hline
\mbox{Group} & & \mbox{Formula} & \\ \hline
G_1 = \langle e \rangle
& N_1(a,b,c) & = & \frac{H(a+b+c)H(a)H(b)H(c)}{H(a+b)H(a+c)H(b+c)} \\
G_2 = \langle \tau \rangle
& N_2(a,a,b) & = & \frac{H_2(2a+b+1)H(a)H_2(b)}{H_2(2a+1)H(a+b)} \\
G_3 = \langle \rho \rangle
& N_3(a,a,a) & = & \frac{H_3(3a+2)H(a)}{H(2a)F_3(3a-2)} \\
G_4 = \langle \tau,\rho \rangle
& N_4(a,a,a)  & = & \frac{H_2(a)H_6(3a+5)}{H_2(2a+1)F_6(3a-2)} \\
G_5 = \langle \kappa \rangle
& N_5(2a,2b,2c) & = & N_1(a,b,c)^2 \\
& N_5(2a,2b,2c+1) & = & N_1(a,b,c)N_1(a,b,c+1) \\
& N_5(2a+1,2b+1,2c) & = & N_1(a+1,b,c)N_1(a,b+1,c) \\
G_6 = \langle \kappa\tau \rangle
& N_6(a,a,2b)& = & \frac{H_2(2b+1)H_2(2b+2a)H(a)}{H(2b+a)H_2(2a)} \\
G_7 = \langle \kappa,\tau \rangle
& N_7(2a,2a,2b) & = & N_1(a,a,b) \\
& N_7(2a+1,2a+1,2b) & = & N_1(a,a+1,b) \\
G_8 = \langle \kappa\tau,\rho \rangle
& N_8(2a,2a,2a) & = & \frac{F_3(3a-2)H_6(6a) H_2(2a)}{H_4(4a+1)H_4(4a)} \\
G_9 = \langle \kappa,\rho \rangle
& N_9(2a,2a,2a)  & \stackrel{!}{=} & \frac{H_3(3a+1)^2 H(a)^2}{H(2a)^2} \\
G_{10} = \langle \kappa,\tau,\rho \rangle
& N_{10}(2a,2a,2a)  & = & \frac{H_3(3a+1)H(a)}{H(2a)} \\[.5em] \hline
\end{array} $$
\caption{\label{formulas} Enumeration formulas for symmetric plane partitions}
\renewcommand{\arraystretch}{1}
\end{table}

In this paper, we will obtain an expression for each of the ten numbers $N_1,
\ldots,N_{10}$ as the determinant or the Pfaffian of a matrix. We will then
row- and column-reduce the matrices to prove the formulas in
Table~\ref{formulas} in cases 1, 3, (part of) 5, and 9.  The author has also
found a way to evaluate the determinants of the matrices for cases 2, 6, 7, and
8, but remains confounded by cases 4 and 10.  The author feels that the
permanent-determinant method is a good first step towards a unified treatment
of all ten cases, but the entirety of such a treatment remains elusive.

The proofs of cases 4 and 10 are more recent than the results of this
paper and will appear in \cite{Stembridge:enumeration} and
\cite{Andrews:tsscpp}.  These three results conclude the problem of
enumeration of plane partitions in different symmetry classes, although
the problem of $q$-enumeration (see section~\ref{qenum}) remains open
in case 4. Table~\ref{credit} gives an account of who first conjectured
and who first proved the formula in each symmetry class.  In
particular, the formula in case 9, whose proof is the main news of this paper, i
Table~\ref{qcredit} gives the same information for $q$-enumeration.
The tables are largely compiled from \cite{Stanley:symmetries}.  They list
the number of the symmetry class as given by Stanley, the name
in the convention of Stembridge, and the generators in the notation
of this paper.

\begin{table}
\begin{center}
\begin{tabular}{|r|c|c||c|c|} \hline
Stanley& Stembridge & Group & First formulator& First enumerator \\ \hline
1 & P    & $\langle e\rangle$              & MacMahon  & MacMahon \\
2 & S    & $\langle\tau\rangle$            & MacMahon  &
 Andrews,Gordon,Macdonald \\
3 & CS   & $\langle\rho\rangle$            & Macdonald & Andrews \\
4 & TS   & $\langle\rho,\tau\rangle$       & Stanley   & Stembridge \\
5 & SC   & $\langle\kappa\rangle$          & Robbins   & Stanley \\
6 & TC   & $\langle\kappa\tau\rangle$      & Proctor   & Proctor \\
7 & SSC  & $\langle\kappa,\tau\rangle$     & Proctor   & Proctor \\
8 & CSTC & $\langle\kappa\tau,\rho\rangle$ & Robbins   & Mills-Robbins-Rumsey \\
9 & CSSC & $\langle\kappa,\rho\rangle$     & Robbins   & Kuperberg \\
10& TSSC & $\langle\kappa,\tau,\rho\rangle$& Robbins   & Andrews \\ \hline
\end{tabular}
\end{center}
\caption{\label{credit} Conjecturers and provers of plane partition formulas}
\end{table}
\begin{table}
\begin{center}
\begin{tabular}{|r|c|c||c|c|} \hline
Stanley&Stembridge&Group&First $q$-formulator&First $q$-enumerator \\ \hline
1 & P   & $\langle e\rangle$    & MacMahon            & MacMahon \\
2 & S   & $\langle\tau\rangle$  & MacMahon            & Andrews,Macdonald \\
2'& S   & $\langle\tau\rangle$  & Bender-Knuth,Gordon & Andrews,Gordon \\
3 & CS  & $\langle\rho\rangle$  & Macdonald           & Mills-Robbins-Rumsey \\
4'& TS  & $\langle\rho,\tau\rangle$ & Andrews,Robbins & still open \\ \hline
\end{tabular}
\end{center}
\caption{\label{qcredit} $q$-Conjecturers and $q$-provers of plane partition
formulas}
\end{table}

Our plan is to combine two old ideas about combinatorial tilings.  The
first idea has probably been known since antiquity, and is illustrated
in \cite{DT:calissons} and \cite{Robbins:story} and discussed in
great generality in \cite{Thurston:conway}. It is that a plane
partition in a box is equivalent to a tiling by certain rhombuses of a
certain hexagon.  To convert a plane partition to a tiling, we view the
partition as a union of unit cubes in an $a \times b \times c$ box, and
we simply ``draw a picture'' of the partition, together with the back
three sides of the box: 
$$\includegraphics{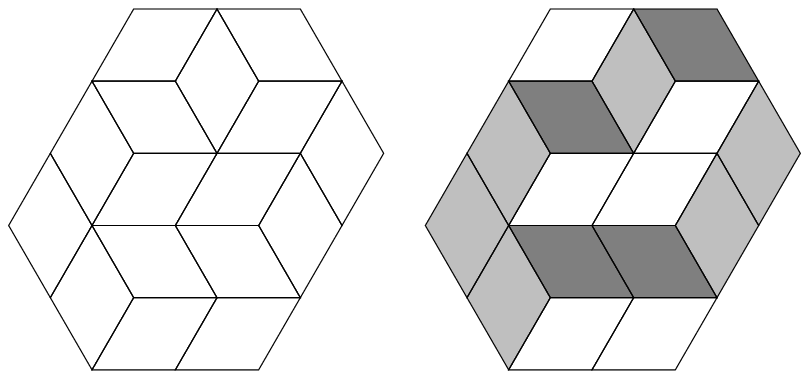}$$
We decree that the hexagon has angles of
$2\pi/3$ and congruent opposite sides of lengths $a$, $b$, and $c$ (we
will call such a hexagon semi-regular), and the rhombuses have angles
of $\pi/3$ and $2\pi/3$ and unit edge length (we will call such a
rhombus a lozenge).  Conversely, we can ``look at'' a tiling and
``see'' the corresponding partition.  We invite the reader to arrive at
this visual proof by staring at some representative rhombus tilings.

The three operations $\tau$, $\rho$, and $\kappa$ have a simple
interpretation in terms of tilings.  The operation $\tau$ is a
reflection about a diagonal between two opposite vertices, $\rho$ is a
rotation by $2\pi/3$ (clockwise by convention), and $\kappa$ is a
rotation by $\pi$.  For example, the number $N_9(2a,2a,2a)$ is also the
number of lozenge tilings of a regular hexagon which are invariant
under rotation by $\pi/3$.

Suppose we are given a tiling of a semi-regular hexagon $H$ by
lozenges.  There is a unique tiling of $H$ by unit equilateral
triangles.  We may form a graph $Z$ whose vertices are these triangles,
with two vertices connected when the triangles are adjacent.  Each
tiling by rhombuses corresponds to a (perfect) matching in $Z$.  This
brings us to the second idea, the permanent-determinant method, which
can be used to count the number of matchings in a bipartite planar
graph, and the Hafnian-Pfaffian method, a generalization for
non-bipartite graphs.  The general method is due to Kasteleyn
\cite{Kasteleyn:crystal} \cite{Kasteleyn:dimer}, but was motivated by
an important special case independently discovered by Kasteleyn
\cite{Kasteleyn:quadratic} and Temperley and Fisher \cite{TF:dimer},
which in turn expanded on ideas of Hurst and Green \cite{HG:ising},
Caianiello \cite{Caianiello:dyson}, and Kac and Ward \cite{KW:ising}.
The permanent-determinant case of the method was noted by Percus 
\cite{Percus:dimer}.

We can use the permanent-determinant method in cases 1, 3, 6, and 8.
The Hafnian-Pfaffian method applies in the remaining cases.

This paper is a result of a collaboration with James Propp, who has maintained
a steady interest in its contents and its completion, and the ideas presented
here grew out of those in \cite{Kuperberg:eklp1,Kuperberg:eklp2}.  In
particular I would like to thank him for correcting many minor errors in the
final version.  I would also like to thank William Jockusch, Richard Stanley,
John Stembridge, and my former advisor, Andrew Casson, for encouragement and
helpful discussions.

\section{\label{hafpfaff} Definition of the general methods}

We begin with a brief overview of four functions of matrices
known as the permanent, determinant, Hafnian, and Pfaffian.

Let $A$ be an $n \times n$ matrix.  The permanent of $A$
is a sum over permutations of $n$ letters:
$$\Per(A) = \sum_{\sigma \in S_n} \prod_{i=1}^n A_{i,\sigma(i)}$$
The determinant is a signed sum of the same type:
$$\Det(A) = \sum_{\sigma \in S_n} (-1)^{\sigma} \prod_{i=1}^n A_{i,\sigma(i)}$$
In general, the permanent of a matrix is hard to compute 
(see \cite{Valiant:permanent}) and satisfies relatively few interesting
identities, while the determinant is a fundamental object of all
mathematics and can be computed in polynomial time.

A (perfect) matching of a set $S$ of $2n$ objects is a set of $n$ disjoint
pairs of elements of $S$.  A matching is ordered or unordered depending
on whether its members are ordered or unordered pairs.  Let $U_n$ be the
set of unordered matchings among $2n$ objects, and let $O_n$ be
the set of ordered matchings.  The Hafnian and the Pfaffian are polynomials
which bear the same relationship to unordered matchings that the permanent
and determinant do to ordered matchings, although the Pfaffian is most
conveniently defined by ordered matchings.  (Both the term ``Hafnian'' and the
polynomial it denotes were devised by Caianiello, who explains that the Hafnian,
like the element Hafnium, is named after the city of Copenhagen.)  Like a
permutation, an ordered matching has a sign:  The matching
$(1,2),(3,4),\ldots,(2n-1,2n)$ is defined to be positive, and if $\sigma$ is a
permutation of $2n$ elements and $m$ is an ordered matching, then the sign of
$\sigma(m)$ is defined to be the product of the signs of $m$ and $\sigma$.

If $A$ is a symmetric $2n \times 2n$ matrix, the Hafnian of $A$ is defined to be
the sum:
$$\Hf(A) = \sum_{m \in U_n} \prod_{\{i,j\} \in m} A_{i,j}$$
If $A$ is an antisymmetric $2n \times 2n$ matrix, the Pfaffian of $A$ is defined to
be the sum:
$$\Pf(A) = {1 \over 2^n} \sum_{m \in O_n} (-1)^m \prod_{(i,j) \in m} A_{i,j}$$
After close scrutiny of this formula, it becomes clear that the use
of oriented matchings is a notational device and the Pfaffian is
really a polynomial with one term (with a leading coefficient of $\pm 1$) 
for each unoriented matching.  The Hafnian is a generalization of the
permanent and the Pfaffian is a generalization of the determinant, because
if
$$A = \left(\begin{array}{c|c} 0 & B \\ \hline B & 0 \end{array} \right)$$
and 
$$A' = \left(\begin{array}{c|c} 0 & B \\ \hline -B & 0 \end{array} \right),$$
then
$$\Hf(A) = \Per(B)$$
and
$$\Pf(A') = \Det(B)$$
Being a generalization of the permanent, the Hafnian is also in general
intractible.  But the Pfaffian satisfies the important identities
$$\Pf(B^T A B) = \Det(B)\Pf(A)$$
and
$$\Pf(A)^2 = \Det(A),$$
which are closely related to each other.  For a proof of these assertions
about the Pfaffian, see \cite{Kasteleyn:crystal}.  It is also an
instructive exercise to develop both the above definition of the Pfaffian
and its properties from the following alternative definition:  If the
anti-symmetric matrix $A$ is viewed as the matrix of a 2-form $\omega_A \in 
\bigwedge^2 \R^{2n}$, and $\mu$ is the standard volume form on $\R^{2n}$,
then:
$$ \omega_A^{\wedge n} = \Pf(A) \mu.$$

Let $Z$ be a bipartite graph with black vertices $b_1,\ldots,b_n$ and white
vertices $w_1\ldots,w_n$.  We define the bipartite adjacency matrix $M_Z$ of $Z$
by setting $(M_Z){i,j}$ to be the number of edges from $b_i$ to $w_j$.  We see
that the permanent of $M_Z$ is simply the number of matchings in $Z$. More
generally, if $Z_w$ is a weighted version of $Z$, we define $(M_{Z_w})_{i,j} $
to be the total weight of all edges from $b_i$ to $w_j$.  The permanent of
$M_{Z_w}$ is the total weight of all matchings in $Z_w$, where the weight of a
matching in $Z_w$ is the product of the weights of the edges of the matching.

Suppose that $Z$ is planar and is given with a specific planar embedding,
and suppose $Z_w$ is a weighted version of $Z$.  We
wish to obtain a new wighting $Z_{w'}$, by changing the sign of some of the
weights,so that:
\bqn \Det(M_{Z_{w'}}) = \pm\Per(M_{Z_w}) \label{perdet} \eqn
on a term-by-term basis.  We define the weights of $Z_w$ to be the 
product of those of $Z_w$ and the weights of a signed graph $Z_{\pm}$.  We say that
$Z_\pm$ is flat if for every face of $Z$ (in the given embedding) that has
$4k$ sides, an odd number of sides are negative in $Z_\pm$; and for
every face that has $4k+2$ sides, an even number of sides negative in $Z_\pm$.

For the Hafnian-Pfaffian case, let $Z$ be a weighted graph which is not
necessarily bipartite.  If $S_Z$ is the usual symmetric, weighted adjacency
matrix of $Z$, then $\Hf(S_Z)$ is the weighted sum of all matchings in $Z$.

To have a Pfaffian, we need an anti-symmetric matrix.  We let $Z_w$ be
an oriented, weighted graph with vertices $v_1,\ldots,v_n$.  We define
the anti-symmetric incidence matrix $A_{Z_w}$ of $Z_w$ by $A_{i,j} =
w_{i,j} - w_{j,i}$, where $w_{i,j}$ is the total weight of all edges
from $v_i$ to $v_j$.

Suppose that $Z$ is an unoriented weighted graph with a given planar
embedding.  We wish to orient the edges of $Z$ to obtain a graph $\vec{Z_w}$
so that:  \bqn \Pf(A_{\vec{Z_w}}) = \pm\Hf(S_{Z_w}) \label{hafpf} \eqn
Once again, we use a consistency rule:  We say that the orientation of
$\vec{Z_w}$ is flat if each face $F$ of $\vec{Z_w}$ has an odd number of
edges which point in the clockwise direction around $F$.

Kasteleyn's theorem states that every planar graph has a flat orientation, and
that equation~\ref{hafpf} holds if $\vec{Z_w}$ is flat.  These results have an
excellent treatment in \cite{Kasteleyn:crystal} which would not be improved if it
were repeated here.  However, just as the determinant is an interesting special
case of the Pfaffian, the permanent-determinant method is an interesting special
case of the Hafnian-Pfaffian method that admits a self-contained treatment that
we present here:

\begin{theorem} Every planar, bipartite graph $Z$ with an even number of vertices
has at least one flat weighting $Z_\pm$. \label{exists} 
\end{theorem}
\begin{proof}
We view $Z$ as embedded in the sphere and we will achieve flatness on the
outside face as well as the others.  (In fact, it is impossible to have
flatness on all but one face on the sphere.)  Let $f_0$ be the number
of faces with $4k$ sides, let $f_2$ be the number of faces with $4k+2$
sides, and let $e$ and $v$ be the number of edges and vertices.  Then $e
= f_2 \bmod 2$ while $v = 0 \bmod 2$.  Therefore the  Euler characteristic
equation $f_0+f_2-e+v = 2$ reduces to $f_2 = 0 \bmod 2$, or in other words, an even number of
faces have $4k+2$ sides.  Choose a matching among these faces and connect
their centers with paths on the sphere which are transverse to $Z$.
If an edge of $Z$ is crossed by $n$ paths, give it a weight of $(-1)^n$
in $Z_\pm$.  Since every path enters and leaves a face with $4k$ sides
the same number of times, $Z_\pm$ is flat at such a face.  For a
face with $4k+2$ sides, the number of entries has the opposite parity
from the number of exits, therefore $Z_\pm$ is flat here also.
\end{proof}

\begin{theorem} If $Z_w$ is an arbitrary weighting of a planar bipartite
graph $Z$, $Z_\pm$ is
a flat signing, and $Z_{w'}$ is the product weighting,
then:
\bq \Det(M_{Z_{w'}}) = \pm\Per(M_{Z_w})\eq
\end{theorem}
\begin{proof} We would like to prove that all terms in $\Det(M_{Z_{w'}})$
have the same relative sign compared to the corresponding terms in
$\Per(M_{Z_w})$.  It suffices to show that all terms in $\Det(M_{Z_\pm})$
have the same sign.  For simplicity we assume that $Z$ has at most one
edge connecting any pair of vertices.

We first note that if a cycle $c$ in $Z$, not necessarily the boundary of a
face, encircles an even number of vertices, then the cycle has an even number
negative edges if it has $4k+2$ sides and an odd number if it has $4k$
sides.  (Proof:  We define $Z_{c,\pm}$ by removing all edges and vertices of
$Z_{\pm}$ outside of $c$.  We know that $Z_{c,\pm}$ is flat at every internal
face and we want to show it is flat on the outside face as well.  To do this we
show that any signed graph $G_\pm$ with an even number of vertices has an even
number of non-flat sides. By the proof of lemma~\ref{exists}, we know that
$G_\pm$ has an even number of faces with $4k$ sides.  If $G_\pm$ has the trivial
signing, these are precisely the sides which are not flat.  But if the sign of an
edge is reversed, the flatness of precisely two sides is reversed.)

Let $m_1$ and $m_2$ be two matchings in $Z$, and let $\sigma_1$ and  $\sigma_2$
be the corresponding permutations.  The set $m_1 \cup m_2 - m_1 \cap
m_2$ is a union of disjoint cycles, each encircling an even number of vertices
of $Z$, and there is a sequence of matchings in $Z$ connecting $m_1$ to $m_2$
such that two consecutive matchings in the sequence differ by only one cycle. We
want to show that the terms for all of these matchings have the same sign. 
Therefore we can assume by induction that $m_1 \cup m_2 - m_1 \cap m_2$ consists
of a single cycle $c$.

Suppose that $c$ has length $2n$.  (By construction it alternates between edges
of $m_1$ and $m_2$ and therefore has even length.)  Then the permutation
$\sigma_2^{-1} \sigma_1$ is an $n$-cycle, which is an even permutation if $n$ is
odd and vice-versa.  Therefore the relative signs of the permutations is exactly
canceled by the relative signs of the edges of $m_1$ and $m_2$, which is given
by the number of negative edges in $c$.
\end{proof}

The rows and columns of the matrices in this paper will in general
be indexed by an unordered index set rather than by the integers
from 1 to $n$ for some $n$.  The Pfaffian of such a matrix is
only defined up to sign, and since the rows and the columns may
have different index sets, the determinant has the same ambiguity.
Therefore we will work instead with the absolute determinant or
absolute Pfaffian (the absolute value of the determinant of Pfaffian),
which are still unambiguous quantities in this case.

\section{Enumeration in the simplest case}

In this section we give a complete proof of the formula for the
number of plane partitions in case 1.  The analysis of case 1 will
set a pattern for the other cases and will also be used as a lemma.

We will explicitly define matrices whose determinants or Pfaffians
are the answers to the ten counting problems.

Let $H(a,b,c)$ denote a semi-regular hexagon with edge lengths
$a$, $b$, and $c$.  Let $T(a,b,c)$ be the tiling of $H(a,b,c)$ by
equilateral triangles of unit edge length.  We can divide the triangles
of $T(a,b,c)$ into two kinds, \N triangles and \D triangles: 
$$\includegraphics{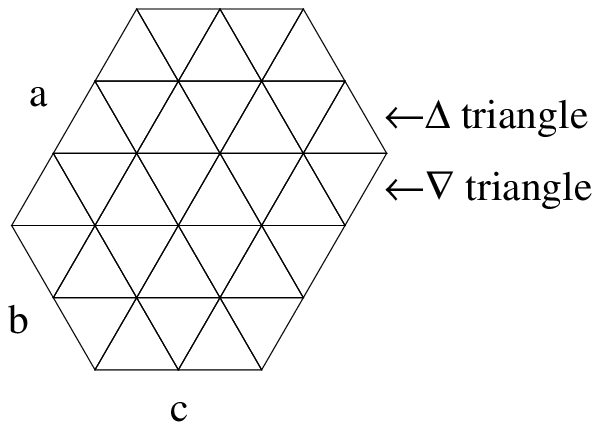}$$
We define a graph $Z(a,b,c)$ with vertex set $T(a,b,c)$, where two
triangles are connected by an edge if and only if they are adjacent.
$Z(a,b,c)$ might look like this:  
$$\includegraphics{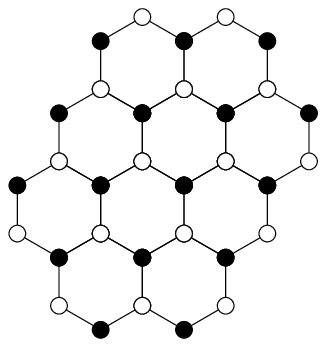}$$
We give $Z(a,b,c)$ the ``default'' weighting, whereby the weight of
each edge is 1.  Evidently, $Z(a,b,c)$ is bipartite, and we can take
the labels \N and \D as a standard bicoloring.

We see that a matching in $Z(a,b,c)$ corresponds to a plane partition
in $B(a,b,c)$.  Thus, $N_1(a,b,c)$ is simply the number of matchings
in $Z(a,b,c)$.  Let $M(a,b,c)$ be the bipartite adjacency matrix of $Z(a,b,c)$.
We decree that the \N triangles are the
rows of $M(a,b,c)$ and the \D triangles are the columns.  To
apply the permanent-determinant method, we need to find a new matrix
whose determinant equals $\Per(M(a,b,c))$.  In this case,
all of the faces of $Z(a,b,c)$ except the outside face have six sides.
Therefore not changing the weights of $Z(a,b,c)$ at all is a flat
sign rule.  In other words,
\bqn |\Det(M(a,b,c))| = \Per(M(a,b,c)) = N_1(a,b,c). \label{perisdet} \eqn


We have obtained a determinant which we must evaluate.  We proceed
by row- and column-reduction and by induction on the area of $H(a,b,c)$.

We embed $H(a,b,c)$ in $H(a+1,b+1,c-1)$ as follows: 
$$\includegraphics{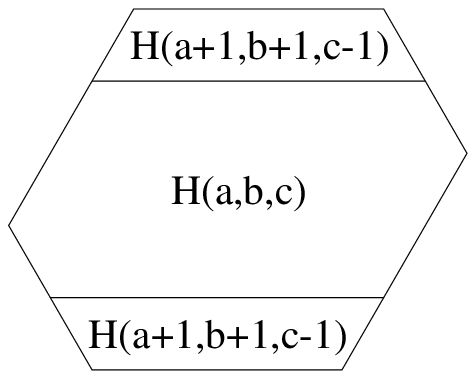}$$
We embed $T(a,b,c)$ in $T(a+1,b+1,c-1)$ and $M(a,b,c)$ in
$M(a+1,b+1,c-1)$ in the corresponding way.  We label the triangles in
and adjacent to
$T(a+1,b+1,c-1)-T(a,b,c)$ with the labels $x_i$, $y_i$, $z_i$, and $l$:  
$$\includegraphics{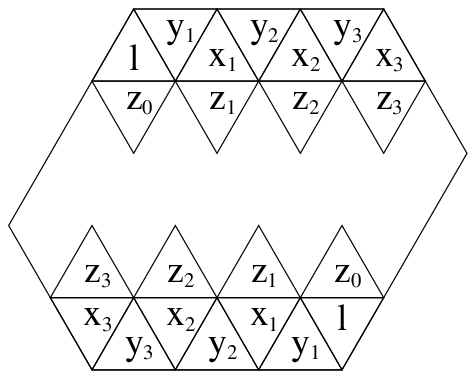}$$
We can use each label twice, once for a \N triangle and once for a \D
triangle.  We use these labels to write out some pieces of the
matrix $M(a+1,b+1,c-1)$.
For example, let $D$ be the submatrix whose rows are
$l$ and the $x_i$'s and whose columns are the $y_i$'s and $z_i$'s.
Note that $D^T$ also appears as a submatrix
in the symmetrical position.  In the case $c=4$, the matrix $D$ looks
like this:
$$\begin{array}{c|rrrrrrr}
    &y_1&y_2&y_3&z_0&z_1&z_2&z_3 \\ \hline
l   & 1 & 0 & 0 & 1 & 0 & 0 & 0 \\
x_1 & 1 & 1 & 0 & 0 & 1 & 0 & 0 \\
x_2 & 0 & 1 & 1 & 0 & 0 & 1 & 0 \\
x_3 & 0 & 0 & 1 & 0 & 0 & 0 & 1 \end{array} $$
If we multiply $D$ by the matrix $E$:
$$\begin{array}{c|rrrr}
    & l &x_1&x_2&x_3 \\ \hline
l   & 1 &-1 & 1 &-1 \\
x_1 & 0 & 1 &-1 & 1 \\
x_2 & 0 & 0 & 1 &-1 \\
x_3 & 0 & 0 & 0 & 1 \end{array}$$
on the left, we obtain the product:
$$\begin{array}{c|rrrrrrr}
    &y_1&y_2&y_3&z_0&z_1&z_2&z_3 \\ \hline
l   & 0 & 0 & 0 & 1 &-1 & 1 &-1 \\
x_1 & 1 & 0 & 0 & 0 & 1 &-1 & 1 \\
x_2 & 0 & 1 & 0 & 0 & 0 & 1 &-1 \\
x_3 & 0 & 0 & 1 & 0 & 0 & 0 & 1 \end{array} $$
Since $E$ has determinant 1, and since the $x_i$ and $l$ rows
and columns are zero outside of $D$ and $D^T$, we can replace
$D$ and $D^T$ by $ED$ and $D^TE^T$ in $M(a+1,b+1,c-1)$ to
obtain a new matrix $M'(a,b,c)$ with the same determinant.

The $x_i$ and $y_i$ rows and columns are irrelevant to the absolute
determinant of the matrix $M'(a,b,c)$, because the $x_i$ row (column)
is the only row (column) in which the $y_i$ column (row) is non-zero,
and the entries at the intersections of these rows and columns are all
$\pm 1$.  If we delete these rows and columns we obtain a matrix
$R(a,b,c)$ which still has the same absolute determinant.  This
matrix is simply $M(a,b,c)$ with an extra row and column labelled $l$.
To complete the description of $R(a,b,c)$, we note that
$R(a,b,c)_{l,z_i} = R(a,b,c)_{z_i,l} = (-1)^i$ and $R(a,b,c)_{l,t} = 0$
or $R(a,b,c)_{t,l} = 0$ where $t$ is any triangle which is not a $z_i$
triangle.


We are left with the problem of reducing the $l$ row.  Since $M(a,b,c)$ is
invertible (because, by induction its determinant is the number of plane
partitions and there is at least one of them), we know that there is some linear
combination of the $M(a,b,c)$ rows which matches the $l$ row exactly except for
the $R(a,b,c)_{l,l}$ entry.  We choose a coefficient $\alpha_t$ for each \N
triangle $t \in T(a,b,c)$ so that:   \bq \sum_t \alpha_t R(a,b,c)_{t,u} +
R(a,b,c)_{l,u} = 0\eq for each \D triangle $u \in T(a,b,c)$.  After adding
$\alpha_t$ times the $t$ row to the $l$ row for each $t$, we obtain a new matrix
$R'(a,b,c)$ for which:  \bq R'(a,b,c)_{l,l} = \sum_t \alpha_t R(a,b,c)_{t,l}\eq
The matrix $R'(a,b,c)$ still has the same absolute determinant as
$M(a+1,b+1,c-1)$ and has the form: $$\left(\begin{array}{c|c} M(a,b,c) &
\begin{array}{c} \cdot \\ \cdot \\ \cdot \\ \cdot \end{array} \\ \hline
\begin{array}{ccccc} 0 & 0 & 0 & \ldots & 0 \end{array} & r \end{array} \right)
$$ where $r = R'(a,b,c)_{l,l}$.  Therefore: \bqn |\Det(M(a+1,b+1,c-1))| =
|R'(a,b,c)_{l,l}||\Det(M(a,b,c))|,\label{rprime}\eqn which will prove the
inductive step.  For convenience, we define $\alpha_i = \alpha_{z_i}$ and
$\beta_i = R(a,b,c)_{l,z_i} = R(a,b,c)_{z_i,l}$, where $i$ is an integer.
Observe that $\beta_i = (-1)^i$.

Let $s$,$t$,$u$, and $v$ be four adjacent triangles in
the middle of $H(a,b,c)$ arranged as follows: 
$$\includegraphics{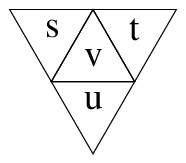}$$
We see that the corresponding $\alpha$ coefficients satisfy a
symmetric Pascal triangle rule:
\bq \alpha_s+\alpha_t+\alpha_u = -R(a,b,c)_{l,v} = 0\eq
The $\beta_i$ coefficients may be incorporated into this rule,
because if $s$,$t$, and $z_i$ are arranged as follows: 
$$\includegraphics{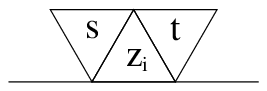}$$
then:
\bq \alpha_s+\alpha_t+R(a,b,c)_{l,z_i} =
\alpha_s+\alpha_t+\beta_i = 0\eq
If we put each coefficient $\alpha_t$ inside the triangle $t$, and
each $\beta_i$ under the \D triangle $z_i$, we see that each number
is the negative of the sum of the two above it: 
$$\includegraphics{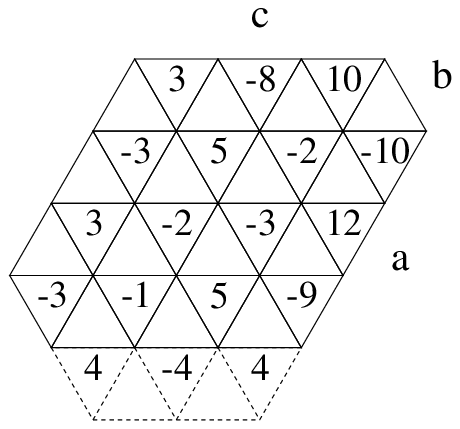}$$
(In this figure the coefficients have all been multiplied by
4 to clear denominators.)
Thus, we get a set of linear equations for the $\alpha_i$'s in terms of
$\beta_i$'s:
\bqn\beta_{c-1-i} = (-1)^{a+b} \sum_j {a+b \choose b+i-j} \alpha_j
\label{case1beta}\eqn
Strictly speaking, the summation over $j$ should be from $0$ to $c-1$
and the equation is only valid for $i$ from $0$ to $c-1$.  However, if
we declare that $\alpha_i = 0$ and $\beta_i$ is an indeterminate for
$i<0$ and $i\ge c$, then the equation will be valid for all $i$ and the
summation can be extended to go from $-\infty$ to $\infty$.  (Recall
that out-of-bounds binomial coefficients are zero by common convention
if the top index is a non-negative integer and the bottom index is an
integer.  Among other things this implies that the $\beta_i$'s must
have finite support if the $\alpha_i$'s do.)  All subsequent summations
will run over all of $\Z$ unless otherwise noted.

We also get an expression for the final
answer in terms of the $\alpha_i$ and $\beta_i$ coefficients:
\bq R'(a,b,c)_{l,l} = \sum_i \alpha_i\beta_i.\eq

The author originally obtained the solution to these equations
with the aid of MACSYMA\TM, but with hindsight we may use
generating functions to make the computation look easy.  We will work in the
formal power series ring $\R[[x,y]][1/x,1/y]$ (which should
not be confused with $\R((x,y))$).  Keeping in mind that the $\alpha_i$
and $\beta_i$ are implicitly functions of $a$, $b$, and $c$, we define
these generating functions:
\begin{eqnarray*}
A(a,b,c;x) & = & \sum_i \alpha_i x^i\\
A(a,b;x,y) & = & \sum_{c=0}^\infty A(a,b,c;x) y^{c-1} \\
M(a,b;x) & = & {(-1)^{a+b}x^{-b}(1+x)^{a+b}} \\
M(a,b;x,y) & = & M(a,b;x) = (-1)^{a+b}x^{-b}y^{-a-b}(xy+y)^{a+b}\\
B(a,b,c;x) & = & \sum_i \beta_{c-1-i} x^i \\
B(a,b;x,y) & = & \sum_{c=0}^\infty B(a,b,c;x) y^{c-1}
\end{eqnarray*}
Then the equation
\bq B(a,b,c;x) = M(a,b;x)A(a,b,c;x),\eq
which is equivalent to
\bq B(a,b;x,y) = M(a,b;x,y)A(a,b;x,y),\eq
is also equivalent to equation~\ref{case1beta}.  By the definition
of the $\beta_i$ coefficients,
\bq B(a,b;x,y) = {1 \over (1-xy)(1+y)} + \ldots \eq
where the ellipses here and below represent indeterminate or irrelevant terms
with either more factors of $x$ than $y$ or negative powers of $x$, and
$A(a,b;x,y)$ is to be determined.  Proceeding directly to the answer, we set:
\bq A(a,b;x,y)={(-1)^{a+b} \over (1-xy)^{a+1}(1+y)^{b+1} {a+b \choose a}}\eq
We apply the binomial expansion theorem:
\begin{eqnarray*}
M(a,b;x,y)A(a,b;x,y) & = &
{((1+y)-(1-xy))^{a+b} \over x^by^{a+b}(1-xy)^{a+1}(1+y)^{b+1}{a+b\choose a}} \\
& = & \sum_{k=0}^{a+b}{{a+b \choose k}(-1)^k(1-xy)^k(1+y)^{a+b-k} \over
x^by^{a+b}(1-xy)^{a+1}(1+y)^{b+1}{a+b \choose a}}
\end{eqnarray*}
A given term in this sum with $k<a$, when expanded in $x$ and $y$, has more
factors of $x$ than $y$, while a term with $k>a$ has only negative
powers of $x$.  Therefore, only the $k=a$ term contributes to the coefficients
of $B(a,b;x,y)$ that we have constrained.  We compute:
\begin{eqnarray*}
M(a,b;x,y)A(a,b;x,y) & = & \ldots+{(-1)^a\over x^by^{a+b}(1-xy)(1+y)}+\ldots \\
& = & \ldots + {1 \over (1-xy)(1+y)} + \ldots \\
& = & B(a,b;x,y),
\end{eqnarray*}
as desired.  In the second line of the equation we again discard
terms with more factors of $x$ than $y$ or negative powers of $x$.

Finally, we recall:
\bq R'(a,b,c)_{l,l} = \sum_i\alpha_i\beta_i = A(a,b,c;-1)\eq
and we compute:
\bq \sum_{c=0}^\infty A(a,b,c;-1) y^{c-1} = A(a,b;-1,y) = {(-1)^{a+b} \over
{a+b \choose a}(1+y)^{a+b+2}}\eq
Using the power series expansion:
\bq {1 \over (1-x)^k} = \sum_n {n+k-1 \choose n} x^n,\eq
we obtain:
\bqn R'(a,b,c)_{l,l} = (-1)^{a+b+c-1}{a+b+c \choose c-1}
\bigg/ {a+b \choose a},\label{case1answer}\eqn
or equivalently, by equations \ref{perisdet} and \ref{rprime}:
\bq {N_1(a+1,b+1,c-1) \over N_1(a,b,c)} = {a+b+c \choose c-1}
\bigg/ {a+b \choose a},\eq
which is the desired ratio.

For comparison with the other cases of enumeration, we introduce an extra
subscript for many of the  quantities defined here.  We will use $\alpha_{1,i}$
for $\alpha_i$, $\beta_{1,i}$ for $\beta_i$, $A_1$ for $A$, and $B_1$ for $B$.

\section{Determinants and Pfaffians for matchings with symmetries}

In this section, we will analyze the set of matchings of a graph which
are invariant under some group action.  We will then apply this analysis to
the particular group actions on the graph $Z(a,b,c)$ under
consideration.  The result for each group will be a planar graph whose
matchings correspond to plane partitions with the equivalent symmetry.
We can then use the permanent-determinant or Hafnian-Pfaffian method to
reduce each of the ten enumerations to questions about matrices.

Suppose that $Z$ is a finite graph with a distinguished vertex called
``bachelorhood''.  We define a matching with bachelors in $Z$ to be a
matching in which every vertex is matched to exactly one of its
neighbors, except for the bachelorhood vertex, which can be matched to
any subset of its neighbors.  (Note that there may be more than one
edge between the bachelorhood vertex and another vertex.  We
define the bachelor's degree of a vertex to be the number of edges
connecting it to the bachelorhood vertex.) The concept of matchings
with bachelors is a generalization of ordinary matchings in the sense
that if $Z$ is an ordinary graph, we let $Z'$ be $Z$ together
with an isolated bachelorhood vertex, and a matching in $Z$ will
correspond to a matching with bachelors in $Z'$.

Given a graph $Z$ with a bachelorhood vertex, we wish to construct a
graph $Z'$ such that ordinary matchings in $Z'$ correspond to matchings
with bachelors in $Z$.  For this purpose we can replace the bachelorhood
vertex of $Z$ by a subgraph which consists of either a row of
triangles: 
$$\includegraphics{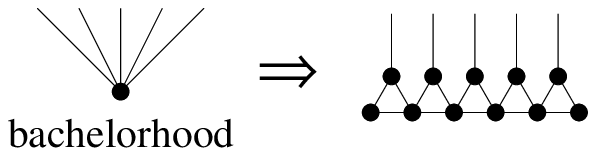}$$
or a row of triangles and an edge: 
$$\includegraphics{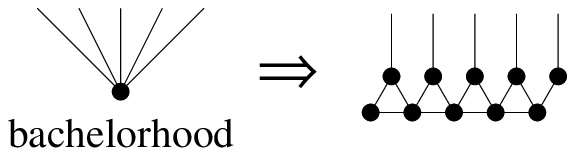}$$
The first subgraph has the property that if an odd number of vertices
are matched to the outside, the remainder can be matched among themselves
in a unique way, and the second subgraph has the analogous property
for an even number of vertices.  (The proof of this fact is left as an
exercise to the reader.)  Thus, we replace the bachelorhood
vertex by one of these two graphs in such a way that the total number
of vertices is even and call the result $Z'$.  The graph $Z'$ has exactly
one matching for each matching with bachelors of $Z$.  Moreover, if $Z$
is a planar graph, it possible to construct $Z'$ to be planar also.

Suppose that $Z$ is a graph with a bachelorhood vertex (which may be
isolated) and a group $G$ acts on $Z$.  That is, each element of $G$ is
some permutation of the vertices that fixes the bachelorhood vertex and
edges of $Z$ that preserves the relation of an edge containing a
vertex.  We wish to construct a modified quotient graph $Z//G$ whose
matchings are bijective with $G$-invariant matchings of $Z$.  If $e$ is
an edge reversed by $G$, i.e. there exists an element of $G$ which
fixes $e$ but switches its endpoints, we replace $e$ by an edge from
each of its endpoints to the bachelorhood vertex.  We perform this
operation for all edges reversed by $G$ and call the result $Z'$.  The
action of $G$ on $Z$ extends to an action on $Z'$, and the
$G$-invariant matchings of $Z'$ can be identified with the
$G$-invariant matchings of $Z$.  Suppose that $e$ is an edge of $Z'$
with a vertex $v$ which is not the bachelorhood vertex such that the
stabilizer of $e$ does not contain the stabilizer of $v$.  Then $e$
cannot be a member of any $G$-invariant matching of $Z'$.  Let $Z''$ be
$Z'$ with all such edges removed.  Let $Z''/G$ be the usual quotient of
a graph by a group action:  The vertices of $Z''/G$ are the $G$-orbits
of vertices in $Z''$, and two vertices are connected by an edge if the
corresponding orbits of vertices are connected by an orbit of edges.
Finally, we define $Z//G$ to be $Z''/G$.  The graph $Z//G$ has the
property that:

\begin{lemma} The matchings of $Z//G$ are in natural bijection with the
$G$-invariant matchings of $Z$.
\end{lemma}

The proof of this lemma is sketched by the definition above.

It is not true in general that if $Z$ is planar then $Z//G$ is also.
Indeed, if $Z$ is the edge graph of a hexagonal prism: 
$$\includegraphics{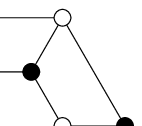}$$ then $Z$
has a fixed-point free involution such that the quotient is $K_{3,3}$,
namely the involution consisting of rotation by 180 degrees and
transposition of the inner and outer hexagons.  However, suppose that
$Z = Z(a,b,c)$ for some $a$, $b$, and $c$, with an isolated
bachelorhood vertex added, and suppose that $G = G_i$ for some $2 \le i
\le 10$.  Then $Z(a,b,c)//G_i$ is indeed a planar graph in these
cases.  We further modify $Z//G$ to remove the bachelorhood vertex and
call the result $Z_i(a,b,c)$.  If $i$ is 3 or 9, or if $i$ is 5 and
$a$, $b$, and $c$ all have the same parity, then $G_i$ acts freely on
$Z(a,b,c)$ and the quotient graph is particularly simple.  When
$\kappa\tau \in G_i$, namely in cases 6,7,8, and 10, the vertices and
edges of $Z(a,b,c)$ fixed by $\kappa\tau$ or its conjugates map to a
collection of disjoint edges in $Z_i(a,b,c)$, which can be ignored
since there is only one matching in a graph which is itself a
matching.  When $\tau \in G_i$, namely in cases 2, 4, 7, and 10,
$Z_i(a,b,c)$ has a non-trivial bachelorhood vertex, so it is necessary
to modify it in the manner described above to obtain a planar graph
without a bachelorhood vertex.  Finally, in cases 3, 6, and 8, the end
result is a bipartite graph, and therefore the permanent-determinant
method is applicable.  In the other six cases, the Hafnian-Pfaffian
method must be used.

\section{Determinants and Pfaffians for $q$-enumeration \label{qenum}}

The previous two sections describe how to construct matrices whose
determinants and Pfaffians enumerate plane partitions in each of the
ten symmetry classes.  The permanent-determinant method applied to
weighted graphs yields matrices whose determinants and Pfaffians
$q$-enumerate plane partitions in various ways.  In particular, the
natural $q$-enumerations described in \cite{Stanley:symmetries} for
cases 1,2,3, and 4 can all be expressed this way.  In this section we
give a determinant for $q$-enumeration in the simplest case, case 1, as
an illustration of the general principle.

We define the $q$-weight of a plane partition to be $q^n$ if the plane
partition has $n$ elements.  We wish to find the total $q$-weight of
all plane partitions in an $a \times b \times c$ box.  We say that two
plane partitions in such a box differ by an elementary move if they
are the same except for one cube.  Observe that we can connect any
two plane partitions in a box by a sequence of elementary moves.
If two plane partitions differ by an elementary move, the corresponding
matchings in $Z(a,b,c)$ differ on only one hexagon: 
$$\includegraphics{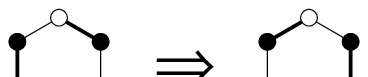}$$
In fact, there are two kinds of elementary moves:  A raising elementary
move is in which a cube is added, while a lowering elementary move
is one in which a cube is deleted.

We wish to weight the edges of $Z(a,b,c)$ in such a way that the weight
of a matching agrees with the weight of the corresponding plane partition.
If we find such a weighting with the property that if we perform a
raising elementary move on a matching, the weight goes up by a factor of $q$,
then the total weight of all matchings will equal the total weight of all
plane partitions up to a constant factor.  That constant factor is the
weight of the matching corresponding to the empty plane partition, because
the $q$-weight of the empty plane partition is 1.  To accomplish this
we give all slanted edges in $Z(a,b,c)$ and all of the left-most vertical
edges a weight of 1.  The weight of any other vertical edge is $q$ times
the weight of the edge immediately to the left of it.  The result is
of the following form: 
$$\includegraphics{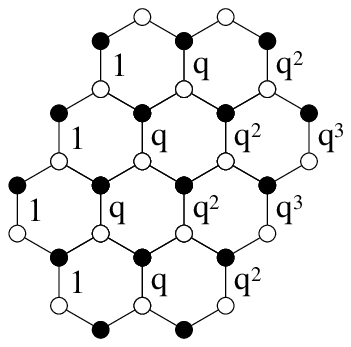}$$
The bipartite adjacency matrix of this graph has the desired determinant
up to a normalization factor.

The total $q$-weight is also an interesting statistic in cases 2, 3,
and 4, and in each of these cases the same method yields a matrix whose
determinant or Pfaffian is the desired statistic.  However, there is no
known formula for the total $q$-weight of all plane partitions in case
4 in this sense.  For the purpose of understanding case 4, there is an
alternative weighting scheme:  The $q$-weight of a plane partition is
defined to be $q^n$ if the plane partition contains $n$ orbits under
the required symmetry group.  The total orbit-$q$-weight satisfies a
nice formula in cases 2 and 4, and again in those two cases as well as
in case 3 the above method yields a determinant or Pfaffian.  The
orbit-$q$-enumeration of plane partitions in case 4 remains the only
open enumeration in the standard conjectures about symmetric plane
partitions.

\section{Enumeration in a previously open case}

The ultimate goal of this section is to prove the formulas for
$N_3(a,a,a)$, $N_5(2a,2b,2c)$, and $N_9(2a,2a,2a)$ by counting
matchings in $Z_3(a,a,a)$, $Z_5(2a,2b,2c)$, and $Z_9(2a,2a,2a)$.  The
hard part of the proof will be row- and column-reduction.  The proof
in case 1 will be used as a template and as a lemma.  In addition, the
three cases treated here together with case 1 form an analogy:
$$ \begin{array}{ccc}
\mbox{case 1} & \mbox{is to} &  \mbox{case 3} \\
\mbox{is to}  & \mbox{as}    &  \mbox{is to} \\
\mbox{case 5} & \mbox{is to} &  \mbox{case 9}
\end{array}$$

Recall that the vertices of $Z(a,b,c)$ are the triangles in a tiling
$T(a,b,c)$ of a hexagon $H(a,b,c)$: 
$$\includegraphics{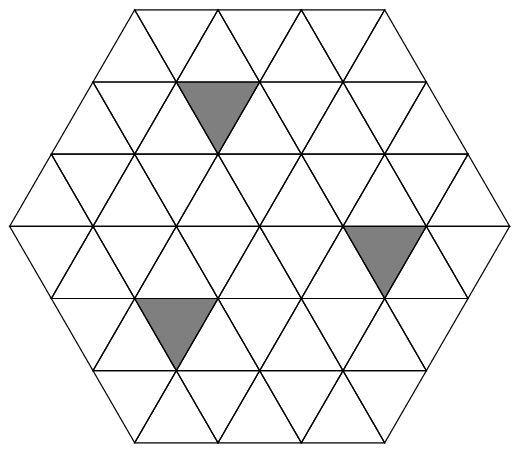}$$
Since $Z_3(a,a,a)$ is a quotient graph of $Z(a,a,a)$, its vertices are
$\rho$-orbits $\{t,\rho(t),\rho^2(t)\}$ of elements of $T(a,a,a)$.
(One of the $\rho$-orbits in the picture is shaded.) Moreover, two
vertices of $Z_3(a,a,a)$ are connected by an edge if and only if they
are adjacent as $\rho$-orbits.  (Note that there is a pair of
$\rho$-orbits in the middle which are connected by two edges.)  Observe
that elements of a $\rho$-orbit are either all \N triangles or \D
triangles, so we can speak of \N $\rho$-orbits and \D $\rho$-orbits.
In other words, the bicoloring of $Z(a,a,a)$ carries over to a
bicoloring of $Z_3(a,a,a)$.  As in case 1, the trivial sign rule is
flat, so if we let $M_3(a,a,a)$ be the bipartite adjacency matrix
of $Z_3(a,a,a)$, its determinant is $N_3(a,a,a)$.  We decree that the
\N $\rho$-orbits should index the rows of $M_3(a,a,a)$ and the \D
$\rho$-orbits should index the columns.

The vertices of $Z_5(2a,2b,2c)$ are pairs of triangles $\{t,\kappa(t)\}$, where
one triangle is a \N triangle and the other is a \D triangle.  Since the graph
$Z_5(2a,2b,2c)$ is not bipartite, we will need to examine its antisymmetric
adjacency matrix under some flat orientation rule, and such a matrix
$M_5(2a,2b,2c)$ would have one row and one column for every vertex of
$Z_5(2a,2b,2c)$.  Instead of having the rows and columns indexed by these
vertices, we can as in case 1 have the rows of $M_5(2a,2b,2c)$ indexed by \N
triangles and the columns by \D triangles, because each \N triangle corresponds
to a unique pair and vice-versa, and the same is true of the \D triangles. With
this indexing, we can take advantage of the fact that $M_5(2a,2b,2c)$ (with any
orientation for $Z_5(2a,2b,2c)$) is simultaneously the bipartite adjacency
matrix for some (non-flat) signed version of $Z(2a,2b,2c)$. The signs must be
chosen so that the an edge from $t$ to $u$ is the negative of the edge from
$\kappa{t}$ to $\kappa{u}$, and it must also correspond to a flat orientation. 
One such weighting is characterized by the following map: 
$$\includegraphics{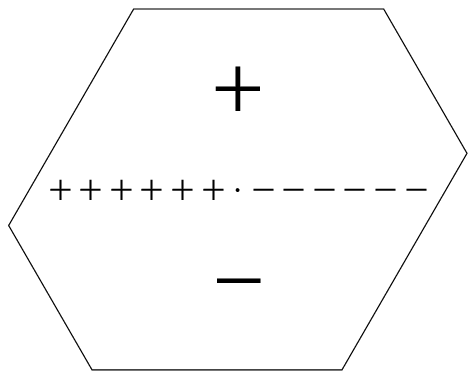}$$ In this
weighting, the edge between two adjacent triangles which are below the center of
$H(2a,2b,2c)$ has weight $-1$, as does the edge between two triangles which are
to the right of the center if one is below the center and one is above the
center.  All other edges have weight 1.  We let $M_5(2a,2b,2c)$ be the bipartite
adjacency matrix of this weighted graph.  It is antisymmetric and its Pfaffian
is $N_5(2a,2b,2c)$.

The same analysis holds for $Z_9(2a,2a,2a)$, except that
the antisymmetric adjacency matrix of an oriented version of this graph
is simultaneously the bipartite adjacency matrix of a weighted version of
$Z_3(2a,2a,2a)$.  We give $Z_3(2a,2a,2a)$ the following weighting: 
$$\includegraphics{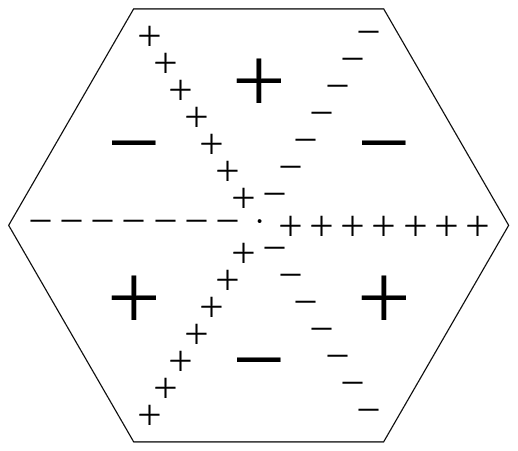}$$
and we let $M_9(2a,2a,2a)$ be the bipartite adjacency matrix of the result.
One way to describe this weighting is to say that if $t$ and $u$ are adjacent
triangles, then:
\bq M_9(2a,2a,2a)_{\{t,\rho(t),\rho^2(t)\},\{u,\rho(u),\rho(u^2)\}}
= M_5(2a,2a,2a)_{t,u} M_5(2a,2a,2a)_{\rho(t),\rho(u)}
M_5(2a,2a,2a)_{\rho^2(t),\rho^2(u)} \eq
The only exception to this rule is in the very center:  If $t$ and $u$ are
the two $\rho$-orbits of triangles in the very center of $H(2a,2a,2a)$,
then $M_9(2a,2a,2a)_{t,u} = 0$.

We embed $M_3(a,a,a)$ in $M_3(a+1,a+1,a+1)$,
$M_5(2a,2b,2c)$ in $M_5(2a+2,2b+2,2c-2)$, and $M_9(2a,2a,2a)$ in
$M_9(2a+2,2a+2,2a+2)$.  The corresponding embeddings of semi-regular
hexagons look like this: 
$$\includegraphics{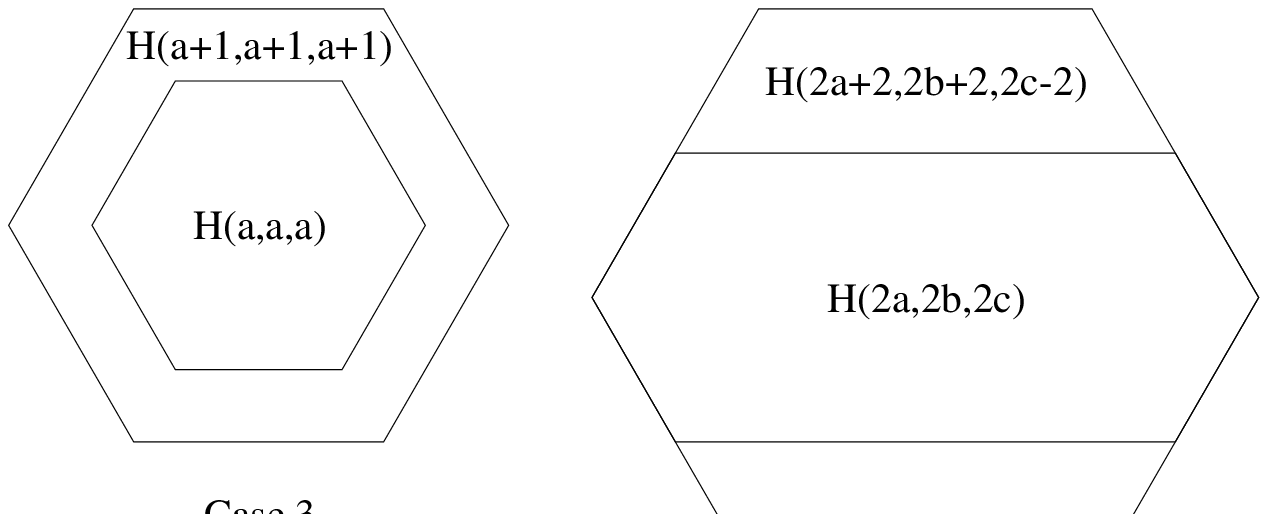}$$
In each case, we label many rows and columns of the larger matrix for
the purpose of row-reduction by labeling the corresponding triangles
or $\rho$-orbits: 
$$\includegraphics{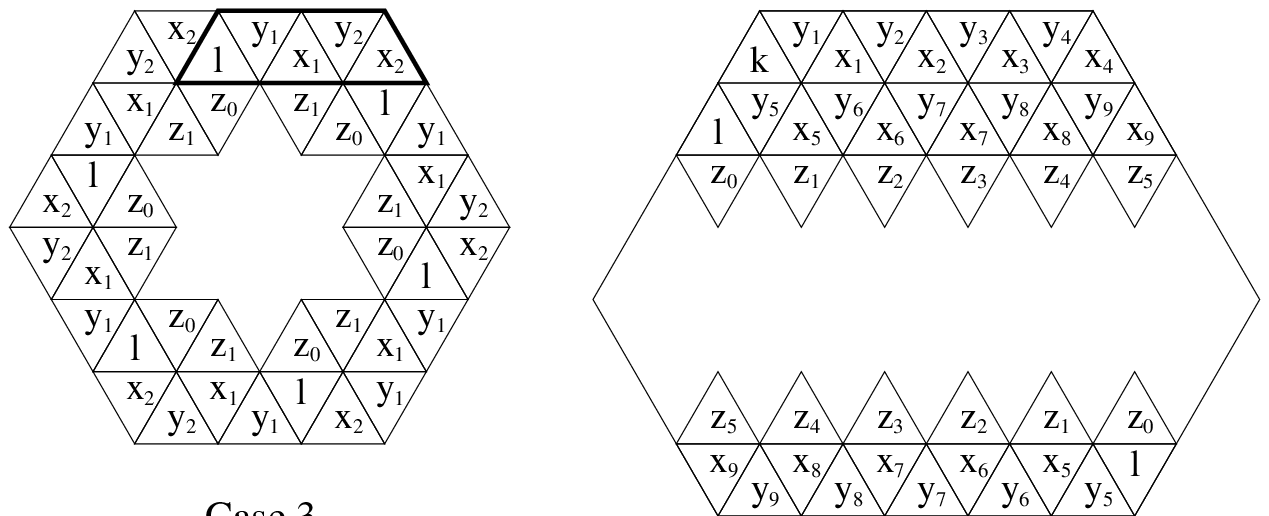}$$

With these labels, the row-reduction algorithm has the same basic plan
as in case 1:   For $i = 3,5,$ or $9$, we define a matrix $R_i(a,b,c)$
by multiplying the $x_i$ rows and columns, the $l$ row and column,
and (if it exists), the $k$ row and column by a certain matrix with
determinant 1.  In the result the $x_j$ and $y_j$ rows and columns
will be irrelevant, and the only non-zero entries of the $l$ and $k$
rows will be in the $z_j$, $l$, and $k$ columns.  However, the actual
values of those entries will be somewhat different.  In all three
cases, we define $\beta_{i,j} = R_i(a,b,c)_{l,z_j}$ and in cases 5
and 9 we also define $\gamma_{i,j} = R_i(a,b,c)_{z_j,k}$.

In the final row-reduction step we will find a linear combination of the
rows of $M_i(a,b,c)$ which matches the $l$ row exactly except for
the $R(a,b,c)_{l,l}$ coefficient and, if it is defined, the $R(a,b,c)_{l,k}$
coefficient.  We choose a coefficient $\alpha_{i,t}$ for each row
$t$ of $M_i(a,b,c)$ so that:
\bq \sum_t \alpha_{i,t} R(a,b,c)_{t,u} + R(a,b,c)_{l,u} = 0 \eq
for each column $u$.  In case 3, we add $\alpha_{3,t}$ times the $t$
row to the $l$ for each $t$ to obtain a matrix $R'_3(a,b,c)$.  Just
as in case 1, the result is:
\bq |\Det(M_3(a+1,a+1,a+1))| = |R'_3(a,a,a)_{l,l}| |\Det(M_3(2a,2a,2a))| \eq
In cases 5 and 9, we add $\alpha_{i,t}$ times the $t$ row to the $l$
row and then add $\alpha_{i,t}$ times the $\kappa(t)$ column to the $l$
column to obtain $R'_5(2a,2b,2c)$ or $R'_9(2a,2a,2a)$.  All operations
performed to obtain these two matrices will be symmetric with respect
to rows and columns, and as a result these two matrices must be
antisymmetric as $M_5(2a,2b,2c)$ and $M_9(2a,2a,2a)$ are.  The matrix
$R'_5(2a,2b,2c)$ will have this form:
$$ \left(\begin{array}{c|c}
M_5(2a,2b,2c) & \begin{array}{cc}
0 & \cdot \\ 0 & \cdot \\ 0 & \cdot \\
\vdots & \vdots \\ 0 & \cdot \end{array} \\ \hline
\begin{array}{ccccc}
0&0&0&\ldots&0 \\
\cdot&\cdot&\cdot&\ldots&\cdot\end{array} &
\begin{array}{cc} 0 & r \\ -r & 0 \end{array} \end{array} \right) $$
where $r = R'_5(2a,2b,2c)_{l,k}$.  From the form of the matrix
we obtain a relation for the absolute Pfaffian:
\bq |\Pf(M_5(2a+2,2b+2,2c-2))| = |R'_5(2a,2b,2c)_{l,k}| |\Pf(M_5(2a,2b,2c))| \eq
and by the same argument for case 9,
\bq |\Pf(M_9(2a+2,2a+2,2a+2))| = |R'_9(2a,2a,2a)_{l,k}| |\Pf(M_9(2a,2a,2a))| \eq

To derive the actual entries of the matrices $R_i(a,b,c)$ and
$R'_i(a,b,c)$ and the coefficients $\alpha_{i,j} = \alpha_{i,z_j}$, we
proceed on a case-by-case basis.  As before, once the $\alpha_{i,j}$'s
are chosen, the values of the other $\alpha_{i,t}$'s can be determined
in a simple manner.

\subsection{Case 3}

In case 3, the submatrix indexed by $x_i$ and $l$ rows and $y_i$
and $z_i$ columns looks the same as in case 1, but we must include
the $l$ and $x_a$ columns because $M_3(a+1,a+1,a+1)_{l,x_a}$ and
$M_3(a+1,a+1,a+1)_{x_a,l}$
are non-zero.  Here is the submatrix $D$ which covers all non-zero
entries of the $l$ and $x_i$ rows in the case $a = 3$:

$$ \begin{array}{c|rrrrrrrr}
    &x_3&y_1&y_2&y_3&z_0&z_1&z_2& l \\ \hline
l   & 1 & 1 & 0 & 0 & 1 & 0 & 0 & 0 \\
x_1 & 0 & 1 & 1 & 0 & 0 & 1 & 0 & 0 \\
x_2 & 0 & 0 & 1 & 1 & 0 & 0 & 1 & 0 \\
x_3 & 0 & 0 & 0 & 1 & 0 & 0 & 0 & 1 \end{array} $$

If we multiply $D$ by the matrix $E$:

$$\begin{array}{c|rrrr}
    & l &x_1&x_2&x_3 \\ \hline
l   & 1 &-1 & 1 &-1 \\
x_1 & 0 & 1 &-1 & 1 \\
x_2 & 0 & 0 & 1 &-1 \\
x_3 & 0 & 0 & 0 & 1 \end{array}$$

on the left, we obtain the product:

$$\begin{array}{c|rrrrrrrr}
    &y_1&y_2&y_3&z_0&z_1&z_2& l  \\ \hline
l   & 0 & 0 & 0 & 1 &-1 & 1 &-1 \\
x_1 & 1 & 0 & 0 & 0 & 1 &-1 & 1 \\
x_2 & 0 & 1 & 0 & 0 & 0 & 1 &-1 \\
x_3 & 0 & 0 & 1 & 0 & 0 & 0 & 1 \end{array} $$

Since $E$ has determinant 1, and since the $x_i$ and $l$ rows
are zero outside of $D$, we can replace
$D$ by $ED$ in $M_3(a+1,a+1,a+1)$ without altering the
determinant.  However, notice that in case 3 there is an intersection
of a single entry, $M_3(a,a,a)_{l,l}$, between $D$ and $D_T$,
and multiplying $D$ on the left by $E$ alters this entry.
Nevertheless we still multiply on the right by $E^T$.  We again
delete the $x_i$ and $y_i$ rows and columns to obtain $R_3(a,a,a)$.
This matrix consists of $M_3(a,a,a)$ plus the non-zero entries
$\beta_{3,i} = R_3(a,a,a)_{l,z_i} = R_3(a,a,a)_{z_i,l} = (-1)^i$ and also
$R_3(a,a,a)_{l,l} = (-1)^a 2$.  We are left with the derivation
of the $\alpha_{3,i}$'s.

We can interpret the $\alpha_{3,t}$ coefficients as a function
on individual triangles in $T(a,a,a)$ rather than $\rho$-orbits
by the equation:
\bq \alpha_{3,t} = \alpha_{3,\{t,\rho(t),\rho^2(t)\}} \eq
Like the $\alpha_{1,t}$'s, the $\alpha_{3,t}$'s satisfy the symmetric
Pascal rule in $H(a,a,a)$, but the boundary conditions are different.
To extend the symmetric Pascal
rule to the boundary, we must place the $\beta_{3,i}$ coefficients
on all of the triangles in the $\rho$-orbits of $x_i$, not just on
the bottom row: 
$$\includegraphics{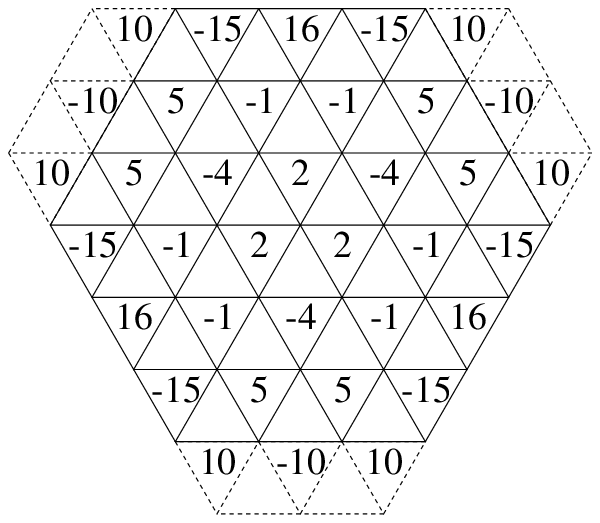}$$
We obtain the solution by setting:
\bq \alpha_{3,t} = \alpha_{1,t} + \alpha_{1,\rho(t)} + \alpha_{1,\rho^2(t)}\eq
for all triangles $t \in T(a,a,a)$.  For integers $i$, we also define
$\alpha'_{1,i} = \alpha_{1,\rho(z_i)}$ and $\alpha''_{1,i} =
\alpha_{1,\rho^2(z_i)}$.

We expand $R'_3(a,a,a)_{l,l}$:
\begin{eqnarray}
R'_3(a,a,a)_{l,l} & = & R_3(a,a,a)_{l,l} + \sum_i (-1)^i \alpha_{1,i} +
\sum_i (-1)^i \alpha'_{1,i} + \sum_i (-1)^i \alpha''_{1,i} \nonumber \\
& = & (-1)^a2 + (-1)^{3a-1} {3a \choose a-1} \bigg/ {2a \choose a} + 2\sum_i
(-1)^i \alpha'_{1,i} \label{case3q}
\end{eqnarray}
by equation~\ref{case1answer} and by the symmetry
$\alpha_{1,\tau\kappa(t)} = (-1)^{a+1}\alpha_{1,t}$, which implies
that $\alpha''_{1,i} = (-1)^{a+1}\alpha'_{1,a-1-i}$.  From the position
of $\rho(z_i)$ in $H(a,a,a)$ (recall that $\rho$ is a clockwise
rotation), we obtain the formula:
\bq \alpha'_{1,i} = (-1)^{a+i}\sum_j {a+i \choose 2a-1-j} \alpha_{1,j}\eq
Therefore:
\bq \sum_i (-1)^i \alpha'_{1,i} = (-1)^a\sum_{i,j} {a+i
\choose 2a-1-j} \alpha_{1,j} = (-1)^a\sum_j {2a \choose 2a-j} \alpha_{1,j}
= (-1)^a\beta_{1,-1}\eq
Recall that $\beta_{1,-1}$ is the $x^a$ coefficient of the polynomial
$B_1(a,a,a;x)$ and the $x^ay^{a-1}$ coefficient of $B_1(a,a;x,y)$.  We
compute:

\bq B_1(a,a,a;x,y) = M_1(a,a;x,y)A_1(a,a,a;x,y) = {((1+y)-(1-xy))^{2a}\over
x^ay^{2a}(1-xy)^{a+1}(1+y)^{a+1}{2a \choose a}}\eq
\bq  = \sum_k{{2a \choose k}(-1)^k(1-x)^k(1+y)^{2a-k} \over
x^ay^{2a}(1-xy)^{a+1}(1+y)^{a+1}{2a \choose a}}\eq
Both the $k=0$ term and the $k=a$ term, but no other terms, contribute to
the $x^ay^{a-1}$ coefficient:
\bq B_1(a,a,a;x,y) = {(1+y)^{a-1} \over x^ay^{2a} (1-xy)^{a+1}{2a\choose a}} +
{(-1)^a \over x^ay^{2a}(1-xy)(1+y)} + \ldots \eq
\bqn  = {x^a y^{a-1}} \left({3a\choose a}\bigg/ {2a\choose a}
-1\right) + \ldots \label{case3b1} \eqn
Substituting equation~\ref{case3b1} in equation~\ref{case3q}, we have:
\bq R'_3(a,a,a)_{l,l} = (-1)^a\biggl(2-{3a\choose a-1}\bigg/{2a\choose a}+
2{3a\choose a}\bigg/{2a\choose a}-2\biggr) = (-1)^a{3a+2\over a}{3a\choose
a-1}\bigg/ {2a\choose a}\eq
This is the desired ratio up to sign.

\subsection{Case 5}

In the two Pfaffian cases, the $x_i$ triangles near a given edge of the hexagon
form two rows rather than one and as a result the initial row-reduction step is
more complicated.  As before, we consider the submatrix $D$ of
$M_5(2a+2,2b+2,2c-2)$ consisting of the $x_i$, $l$, and $k$ rows and the $y_i$
and $z_i$ columns.  Here is $D$ when $c=2$:

$$ \begin{array}{c|rrrrrrrrrrrrrrrr}
    &y_1&y_2&y_3&y_4&y_5&y_6&y_7&y_8&y_9&z_0&z_1&z_2&z_3&z_4&z_5 \\ \hline
k   &-1 & 0 & 0 & 0 &-1 & 0 & 0 & 0 & 0 & 0 & 0 & 0 & 0 & 0 & 0 \\
x_1 &-1 &-1 & 0 & 0 & 0 &-1 & 0 & 0 & 0 & 0 & 0 & 0 & 0 & 0 & 0 \\
x_2 & 0 &-1 &-1 & 0 & 0 & 0 &-1 & 0 & 0 & 0 & 0 & 0 & 0 & 0 & 0 \\
x_3 & 0 & 0 &-1 &-1 & 0 & 0 & 0 &-1 & 0 & 0 & 0 & 0 & 0 & 0 & 0 \\
x_4 & 0 & 0 & 0 &-1 & 0 & 0 & 0 & 0 &-1 & 0 & 0 & 0 & 0 & 0 & 0 \\
l   & 0 & 0 & 0 & 0 &-1 & 0 & 0 & 0 & 0 &-1 & 0 & 0 & 0 & 0 & 0 \\
x_5 & 0 & 0 & 0 & 0 &-1 &-1 & 0 & 0 & 0 & 0 &-1 & 0 & 0 & 0 & 0 \\
x_6 & 0 & 0 & 0 & 0 & 0 &-1 &-1 & 0 & 0 & 0 & 0 &-1 & 0 & 0 & 0 \\
x_7 & 0 & 0 & 0 & 0 & 0 & 0 &-1 &-1 & 0 & 0 & 0 & 0 &-1 & 0 & 0 \\
x_8 & 0 & 0 & 0 & 0 & 0 & 0 & 0 &-1 &-1 & 0 & 0 & 0 & 0 &-1 & 0 \\
x_9 & 0 & 0 & 0 & 0 & 0 & 0 & 0 & 0 &-1 & 0 & 0 & 0 & 0 & 0 &-1
\end{array}$$

Note that $-D^T$ also appears as a submatrix of $M_5(2a+2,2b+2,2c-2)$.
If we multiply $D$ by a matrix $E$ of the form:

$$ \begin{array}{c|rrrrrrrrrrrr}
    & k &x_1&x_2&x_3&x_4&l  &x_5&x_6&x_7&x_8&x_9 \\ \hline
k   & 1 &-1 & 1 &-1 & 1 & 0 &-1 & 2 &-3 & 4 &-5 \\
x_1 & 0 & 1 &-1 & 1 &-1 & 0 & 0 &-1 & 2 &-3 & 4 \\
x_2 & 0 & 0 & 1 &-1 & 1 & 0 & 0 & 0 &-1 & 2 &-3 \\
x_3 & 0 & 0 & 0 & 1 &-1 & 0 & 0 & 0 & 0 &-1 & 2 \\
x_4 & 0 & 0 & 0 & 0 & 1 & 0 & 0 & 0 & 0 & 0 &-1 \\
l   & 0 & 0 & 0 & 0 & 0 & 1 &-1 & 1 &-1 & 1 &-1 \\
x_5 & 0 & 0 & 0 & 0 & 0 & 0 & 1 &-1 & 1 &-1 & 1 \\
x_6 & 0 & 0 & 0 & 0 & 0 & 0 & 0 & 1 &-1 & 1 &-1 \\
x_7 & 0 & 0 & 0 & 0 & 0 & 0 & 0 & 0 & 1 &-1 & 1 \\
x_8 & 0 & 0 & 0 & 0 & 0 & 0 & 0 & 0 & 0 & 1 &-1 \\
x_9 & 0 & 0 & 0 & 0 & 0 & 0 & 0 & 0 & 0 & 0 & 1
\end{array}$$

We obtain:

$$ \begin{array}{c|rrrrrrrrrrrrrrrr}
    &y_1&y_2&y_3&y_4&y_5&y_6&y_7&y_8&y_9&z_0&z_1&z_2&z_3&z_4&z_5 \\ \hline
k   & 0 & 0 & 0 & 0 & 0 & 0 & 0 & 0 & 0 & 0 & 1 &-2 & 3 &-4 & 5 \\
x_1 &-1 & 0 & 0 & 0 & 0 & 0 & 0 & 0 & 0 & 0 & 0 & 1 &-2 & 3 &-4 \\
x_2 & 0 &-1 & 0 & 0 & 0 & 0 & 0 & 0 & 0 & 0 & 0 & 0 & 1 &-2 & 3 \\
x_3 & 0 & 0 &-1 & 0 & 0 & 0 & 0 & 0 & 0 & 0 & 0 & 0 & 0 & 1 &-2 \\
x_4 & 0 & 0 & 0 &-1 & 0 & 0 & 0 & 0 & 0 & 0 & 0 & 0 & 0 & 0 & 1 \\
l   & 0 & 0 & 0 & 0 & 0 & 0 & 0 & 0 & 0 &-1 & 1 &-1 & 1 &-1 & 1 \\
x_5 & 0 & 0 & 0 & 0 &-1 & 0 & 0 & 0 & 0 & 0 &-1 & 1 &-1 & 1 &-1 \\
x_6 & 0 & 0 & 0 & 0 & 0 &-1 & 0 & 0 & 0 & 0 & 0 &-1 & 1 &-1 & 1 \\
x_7 & 0 & 0 & 0 & 0 & 0 & 0 &-1 & 0 & 0 & 0 & 0 & 0 &-1 & 1 &-1 \\
x_8 & 0 & 0 & 0 & 0 & 0 & 0 & 0 &-1 & 0 & 0 & 0 & 0 & 0 &-1 & 1 \\
x_9 & 0 & 0 & 0 & 0 & 0 & 0 & 0 & 0 &-1 & 0 & 0 & 0 & 0 & 0 &-1
\end{array}$$

We replace $D$ by $ED$ and $-D^T$ by $-D^TE^T$ in $M_5(2a+2,2b+2,2c-2)$,
and then delete the $x_i$ and $y_i$ rows and columns to obtain
$R_5(2a,2b,2c)$, which still has the same absolute Pfaffian.
In addition to $M_5(2a,2b,2c)$, the matrix $R_5(2a,2b,2c)$ has
the non-zero entries $\beta_{5,i} = R_5(2a,2b,2c)_{l,z_i}
= -R_5(2a,2b,2c)_{z_i,l} = -(-1)^i$ and $\gamma_{5,i} = -
R_5(2a,2b,2c)_{k,z_i} = i(-1)^i$.

The $\alpha_{5,t}$ coefficients, as a function on the \N triangles of
$T(2a,2b,2c)$, satisfy the symmetric Pascal rule everywhere except
along a ``branch cut'' going from the center of $T(2a,2b,2c)$ to the
right edge, where they satisfy the ordinary Pascal rule: 
$$\includegraphics{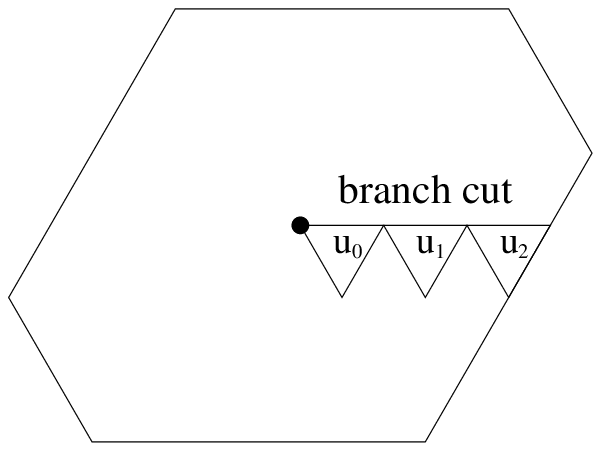}$$
We label the row of triangles just below the branch cut as indicated
and we define $\O{\alpha}_{5,i} = \alpha_{5,u_i}$.

We define the following generating functions:
\begin{eqnarray*}
A_5(2a,2b,2c;x) & = & \sum_i x^i\alpha_{5,i} \\
\O{A}_5(2a,2b,2c;x) & = & \sum_i x^i\O{\alpha}_{5,i} \\
B_5(2a,2b,2c;x) & = & \sum_i x^i\beta_{5,2c-i} \\
\end{eqnarray*}

For convenience in dealing with the branch cut, we define
the function $\theta(n)$ to be $-1$ when $n\ge 0$ and 1 when $n< 0$.
We define the functions $\sigma$ and $\T{\sigma}$ on Laurent
polynomials by:
\begin{eqnarray*}
\sigma(P(x)) & = & \sum_i \theta(i)p_ix^i \\
\T{\sigma}(P(x)) & = & \sum_{i\ne 0} \theta(i)p_ix^i
\end{eqnarray*}
where:
\bq P(x) = \sum_i x^ip_i.\eq

We examine the equations that the $\O{\alpha}_{5,i}$ coefficients satisfy:
\begin{eqnarray}
\O{\alpha}_{5,i} & = & (-1)^{a+b}\theta(i)\sum_j
{a+b \choose a+c+i-j}\alpha_{5,j} \label{alphacase5i}\\
\beta_{5,2c-i} & = & (-1)^{a+b}\sum_j{a+b \choose a-c+i-j} \O{\alpha}_{5,j}
\label{betacase5i}
\end{eqnarray}
We can put equations~\ref{alphacase5i} and \ref{betacase5i} in generating
function form:
\begin{eqnarray*}
\O{A}_5(2a,2b,2c;x) & = & \sigma(x^{-c}M_1(a,b;x)A_5(2a,2b,2c;x)) \\
B_5(2a,2b,2c;x) & = & x^cM_1(a,b;x)\O{A}_5(2a,2b,2c;x)
\end{eqnarray*}
As before, for $i>2c-1$ or $i<0$, $\alpha_{5,i} =0$, but $\beta_{5,i}$ is
indeterminate.

As usual, we jump to the solution (originally found with the aid
of MATHEMATICA\TM) by setting
\bq A_5(2a,2b,2c;x) = (1+x)A_1(a,b,c;x)^2.\eq
Before demonstrating that this is the solution, we define the polynomials
$B_L(a,b,c;x)$ and $B_H(a,b,c;x)$ by:
\bqn B_1(a,b,c;x) = M_1(a,b;x)A_1(a,b,c;x) = (-1)^{c-1}{1 \over 1+x} +
B_L(a,b,c;x) + B_H(a,b,c;x), \label{case5ibh}\eqn
where $B_H(a,b,c;x)$ has no terms with exponents less than $c$, and
$B_L(a,b,c;x)$ consists only of terms with negative exponents.  By
analogy with $B_1(a,b,c;x)$, we also define:
\begin{eqnarray*}
B_H(a,b,c;x) & = & \sum_i \beta_{H,c-1-i} x^i \\
B_L(a,b,c;x) & = & \sum_i \beta_{L,c-1-i} x^i
\end{eqnarray*}
Using the fact that the terms of $A_1(a,b,c;x)$ have degrees from 0 to
$c-1$ inclusive, we can expand $\O{A}_5(2a,2b,2c)$:
\begin{eqnarray*}
\lefteqn{\O{A}_5(2a,2b,2c;x)} \\
& = & \sigma(x^{-c}M_1(a,b;x)(1+x)A_1(a,b,c;x)^2) \\
& = & \sigma(x^{-c}(1+x)A_1(a,b,c;x)({(-1)^{c-1}\over (1+x)}+
B_L(a,b,c;x)+B_H(a,b,c;x))) \\
& = & x^{-c}(1+x)A_1(a,b,c;x)({(-1)^{c-1}\over (1+x)}+
B_L(a,b,c;x)-B_H(a,b,c;x))
\end{eqnarray*}
The result is a cancellation of the ``error terms'' $B_H$ and $B_L$:
\begin{eqnarray}
\lefteqn{x^cM_1(a,b;x)\O{A}_5(2a,2b,2c;x)} \nonumber \\
& = & (1+x)({(-1)^{c-1} \over 1+x} + B_L(a,b,c;x)+B_H(a,b,c;x))
({(-1)^{c-1} \over 1+x} + B_L(a,b,c;x)-B_H(a,b,c;x)) \nonumber \\
& = & {1 \over 1+x} + (-1)^{c-1}2B_L(a,b,c;x) +
(1+x)(B_L(a,b,c;x)^2-B_H(a,b,c;x)^2) \nonumber \\
& = & B_5(2a,2b,2c;x) \label{case5ib5}
\end{eqnarray}
Finally, since $\gamma_{5,i} = i(-1)^i$, we see that:
\bqn R'_5(2a,2b,2c)_{l,k} = \left. {dA_5(2a,2b,2c;x)\over dx}\right|_{x=-1}
= A(a,b,c;-1)^2 = R'(a,b,c)_{l,l}^2 \label{case5ianswer}\eqn

\subsection{Case 9}

Case 9 is similar to case 3 in the sense that the positions in which the
$k$, $l$, and $x_i$ rows are non-zero in the matrix $M_9(2a+2,2a+2,2a+2)$
include the $l$ and $k$ columns as well as the $y_i$ and $z_i$ columns.  We
let $D$ be the submatrix consisting of all of these rows and columns.  Here
is $D$ when $a=4$:

$$ \begin{array}{c|rrrrrrrrrrrrrrrr}
    &y_1&y_2&y_3&y_4&y_5&y_6&y_7&y_8&y_9&z_0&z_1&z_2&z_3& l & k \\ \hline
k   &-1 & 0 & 0 & 0 &-1 & 0 & 0 & 0 & 0 & 0 & 0 & 0 & 0 & 0 & 0 \\
x_1 &-1 &-1 & 0 & 0 & 0 &-1 & 0 & 0 & 0 & 0 & 0 & 0 & 0 & 0 & 0 \\
x_2 & 0 &-1 &-1 & 0 & 0 & 0 &-1 & 0 & 0 & 0 & 0 & 0 & 0 & 0 & 0 \\
x_3 & 0 & 0 &-1 &-1 & 0 & 0 & 0 &-1 & 0 & 0 & 0 & 0 & 0 & 0 & 0 \\
x_4 & 0 & 0 & 0 & 1 & 0 & 0 & 0 & 0 & 1 & 0 & 0 & 0 & 0 & 0 & 0 \\
l   & 0 & 0 & 0 & 0 &-1 & 0 & 0 & 0 & 0 &-1 & 0 & 0 & 0 & 0 & 0 \\
x_5 & 0 & 0 & 0 & 0 &-1 &-1 & 0 & 0 & 0 & 0 &-1 & 0 & 0 & 0 & 0 \\
x_6 & 0 & 0 & 0 & 0 & 0 &-1 &-1 & 0 & 0 & 0 & 0 &-1 & 0 & 0 & 0 \\
x_7 & 0 & 0 & 0 & 0 & 0 & 0 &-1 &-1 & 0 & 0 & 0 & 0 &-1 & 0 & 0 \\
x_8 & 0 & 0 & 0 & 0 & 0 & 0 & 0 & 1 & 1 & 0 & 0 & 0 & 0 & 1 & 0 \\
x_9 & 0 & 0 & 0 & 0 & 0 & 0 & 0 & 0 & 1 & 0 & 0 & 0 & 0 & 0 & 1
\end{array}$$

Some of the entries are positive, because unlike case 5, in case 9 the \N
triangles $x_a$, $x_{2a}$, and $x_{2a+1}$ lie outside of the negative region
in the sign rule map for the matrix $M_9(2a+2,2a+2,2a+2)$.

Note that $-D^T$ also appears as a submatrix of $M_9(2a+2,2a+2,2a+2)$
but that the two matrices intersect in the $2 \times 2$ matrix formed
by $l$ and $k$. We multiply $D$ by a matrix $E$ of the form:

$$ \begin{array}{c|rrrrrrrrrrrr}
    & k &x_1&x_2&x_3&x_4&l  &x_5&x_6&x_7&x_8&x_9 \\ \hline
k   & 1 &-1 & 1 &-1 &-1 & 0 &-1 & 2 &-3 &-4 & 5 \\
x_1 & 0 & 1 &-1 & 1 & 1 & 0 & 0 &-1 & 2 & 3 &-4 \\
x_2 & 0 & 0 & 1 &-1 &-1 & 0 & 0 & 0 &-1 &-2 & 3 \\
x_3 & 0 & 0 & 0 & 1 & 1 & 0 & 0 & 0 & 0 & 1 &-2 \\
x_4 & 0 & 0 & 0 & 0 &-1 & 0 & 0 & 0 & 0 & 0 & 1 \\
l   & 0 & 0 & 0 & 0 & 0 & 1 &-1 & 1 &-1 &-1 & 1 \\
x_5 & 0 & 0 & 0 & 0 & 0 & 0 & 1 &-1 & 1 & 1 &-1 \\
x_6 & 0 & 0 & 0 & 0 & 0 & 0 & 0 & 1 &-1 &-1 & 1 \\
x_7 & 0 & 0 & 0 & 0 & 0 & 0 & 0 & 0 & 1 & 1 &-1 \\
x_8 & 0 & 0 & 0 & 0 & 0 & 0 & 0 & 0 & 0 &-1 & 1 \\
x_9 & 0 & 0 & 0 & 0 & 0 & 0 & 0 & 0 & 0 & 0 &-1
\end{array}$$

This matrix is the same as the one in case 5, except that the $x_a$,
$x_{2a}$, and $x_{2a+1}$ columns have been negated.  The product in this
case is:

$$ \begin{array}{c|rrrrrrrrrrrrrrrr}
    &y_1&y_2&y_3&y_4&y_5&y_6&y_7&y_8&y_9&z_0&z_1&z_2&z_3& l & k \\ \hline
k   & 0 & 0 & 0 & 0 & 0 & 0 & 0 & 0 & 0 & 0 & 1 &-2 & 3 &-4 & 5 \\
x_1 &-1 & 0 & 0 & 0 & 0 & 0 & 0 & 0 & 0 & 0 & 0 & 1 &-2 & 3 &-4 \\
x_2 & 0 &-1 & 0 & 0 & 0 & 0 & 0 & 0 & 0 & 0 & 0 & 0 & 1 &-2 & 3 \\
x_3 & 0 & 0 &-1 & 0 & 0 & 0 & 0 & 0 & 0 & 0 & 0 & 0 & 0 & 1 &-2 \\
x_4 & 0 & 0 & 0 & 1 & 0 & 0 & 0 & 0 & 0 & 0 & 0 & 0 & 0 & 0 & 1 \\
l   & 0 & 0 & 0 & 0 & 0 & 0 & 0 & 0 & 0 &-1 & 1 &-1 & 1 &-1 & 1 \\
x_5 & 0 & 0 & 0 & 0 &-1 & 0 & 0 & 0 & 0 & 0 &-1 & 1 &-1 & 1 &-1 \\
x_6 & 0 & 0 & 0 & 0 & 0 &-1 & 0 & 0 & 0 & 0 & 0 &-1 & 1 &-1 & 1 \\
x_7 & 0 & 0 & 0 & 0 & 0 & 0 &-1 & 0 & 0 & 0 & 0 & 0 &-1 & 1 &-1 \\
x_8 & 0 & 0 & 0 & 0 & 0 & 0 & 0 &-1 & 0 & 0 & 0 & 0 & 0 &-1 & 1 \\
x_9 & 0 & 0 & 0 & 0 & 0 & 0 & 0 & 0 &-1 & 0 & 0 & 0 & 0 & 0 &-1
\end{array}$$

In the matrix $M_9(2a+2,2a+2,2a+2)$, we multiply $D$ by $E$ on the
left, which affects the $l$-$k$ submatrix of $-D^T$, and then we
multiply the modified version of $-D^T$ by $E^T$ on the right.
After multiplying on the left, the $l$-$k$ submatrix is in general:

$$\begin{array}{c|rr}  &k&l\\ \hline k&2a+1 & -2a \\ l & 1 & -1\end{array}$$

and after multiplying on the right, it becomes:

$$\begin{array}{c|rr}  &k&l\\ \hline k& 0 & -2a-1 \\ l & 2a+1 & 0\end{array}$$

We delete the $x_i$ and $y_i$ rows and columns and call the result
$R_9(2a,2a,2a)$.  In addition to the above submatrix and the
submatrix $M_9(2a,2a,2a)$, $R_9(2a,2a,2a)$ also has the non-zero
entries $\beta_{9,i} = R_9(2a,2b,2c)_{l,z_i}
= -R_9(2a,2b,2c)_{z_i,l} = -(-1)^i$ and $\gamma_{9,i} = -
R_9(2a,2b,2c)_{k,z_i} = i(-1)^i$.

As in case 3, we interpret the $\alpha_{9,t}$ coefficients as a function
on triangles in $T(2a,2a,2a)$ by the following equation:
\bq \alpha_{9,t} = \alpha_{9_,\{t,\rho(t),\rho^2(t)\}} \eq
The resulting
coefficients satisfy the symmetric Pascal rule everywhere except along
three branch cuts instead of one: 
$$\includegraphics{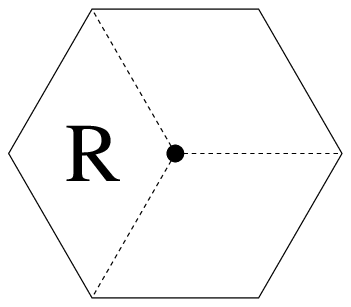}$$ We can recycle the computations
of case 5 by setting $\T{\alpha}_{5,t} = \alpha_{5,t}$
for all $t \in T(2a,2a,2a)$ except those which lie in the region $R$
between the two branch cuts on the left side.  For $t \in R$, we
set $\T{\alpha}_{5,t} = -\alpha_{5,t}$.  We see that the
$\alpha_{9,t}$'s satisfy the same rule as the $\T{\alpha}_{5,t}$'s,
except with different boundary conditions.  Following case 3, we can
obtain the former from the latter by symmetrization:
\bq \alpha_{9,t} = \T{\alpha}_{5,t} +\T{\alpha}_{5,\rho(t)}
+ \T{\alpha}_{5,\rho(t)}\eq
We define:
\bq\alpha'_{5,i} =\T{\alpha}_{5,\rho(z_i)} = \alpha_{5,\rho(z_i)}\eq
\bq\alpha''_{5,i} =\T{\alpha}_{5,\rho^2(z_i)} = -\alpha_{5,\rho^2(z_i)}\eq
(Recall that $\rho$ is a clockwise rotation.)  We expand
$R'_9(2a,2a,2a)_{l,k}$:
\begin{eqnarray}
R'_9(2a,2a,2a)_{l,k} & = & R_9(2a,2a,2a)_{l,k} +
\sum_i i(-1)^i \alpha_{5,i} + \sum_i i(-1)^i \alpha'_{5,i} + \sum_i i(-1)^i
\alpha''_{5,i} \nonumber \\
& = & 2a+1 + {3a \choose a-1}^2 \bigg/ {2a \choose a}^2
+ \sum_i (-1)^i (2i-2a-1) \alpha'_{5,i} \label{case9q}
\end{eqnarray}
from equation~\ref{case5ianswer} and from the symmetry $\alpha_{5,t} =
-\alpha_{5,\kappa\tau(t)}$, which holds for triangles $t$ below the center
of $T(2a,2a,2a)$ and which implies that $\alpha''_{5,i} =
-\alpha'_{5,2a-1-i}$.   From the position of $\rho(z_i)$ in $T(2a,2a,2a)$, we
learn that:
\bq \alpha'_{5,i} = (-1)^i \sum_j{2a+i \choose 4a-1-j} \alpha_{5,j} \eq
Therefore:
\begin{eqnarray}
\sum_i (-1)^i(2i-2a-1)\alpha'_{5,i} & = & \sum_{i,j} (2i-2a-1) {2a+i \choose
4a+1-j} \alpha_{5,j} \nonumber \\
& = & \sum_j \biggl((2a-1){4a \choose 4a-j} - 2{4a \choose
4a+1-j}\biggr) \alpha_{5,j} \nonumber \\
& = & 2\beta_{5,-2} - (2a-1)\beta_{5,-1}\label{case9alpha}
\end{eqnarray}
We want to know $\beta_{5,-1}$ and $\beta_{5,-2}$, which are the $x^{2a}$
and $x^{2a+1}$ coefficients of $B_5(2a,2a,2a;x)$.  If we combine
equations~\ref{case5ibh} and \ref{case5ib5}, we obtain:
\begin{eqnarray}
B_5(2a,2a,2a;x) & = & {1 \over 1+x} - (1+x)B_H(a,a,a;x)^2 + \ldots \nonumber \\
& = & \ldots+x^{2a}-x^{2a+1}+\ldots-
(1+x)(\beta_{H,-1}x^a + \beta_{H,-2}x^{a+1}+\ldots)^2 \nonumber \\
& = & (1-\beta_{H,-1}^2)x^{2a} +
(-1-\beta_{H,-1}(2\beta_{H,-2}+\beta_{H,-1}))x^{2a+1} + \ldots\label{case9b5}
\end{eqnarray}
Putting together equations~\ref{case9q}, \ref{case9alpha}, and \ref{case9b5},
we obtain:
\bq R'_9(2a,2a,2a)_{l,k} = {a^2 \over (3a+1)^2} {3a+1 \choose a}^2 \bigg/
{2a \choose a}^2 - \beta_{H,-1}(4\beta_{H,-2}-(2a-3)\beta_{H,-1}) \eq
We know from equation~\ref{case3b1} that:
\bq \beta_{H,-1} = {3a \choose a} \bigg/ {2a \choose a}
= {2a+1 \over 3a+1} {3a+1 \choose a} \bigg/ {2a \choose a} \eq
We can obtain $\beta_{H,-2}$ in the same way we derived $\beta_{H,-1}$:
\begin{eqnarray*}
B_1(a,a,a;x,y) & = & {(1+y)^{a-1} \over x^ay^{2a}(1-xy)^{a+1} {2a \choose a}}
- {2a(1+y)^a-2 \over x^ay^{2a}(1-xy)^a{2a \choose a}} + \ldots +
{(-1)^a \over x^ay^{2a}(1+y)(1-xy)} + \ldots \\
& = & (a-1)x^{a+1}y^{a-1}{3a+1 \choose a} \bigg/ {2a \choose a}
- 2ax^{a+1}y^{a-1}{3a \choose a-1} \bigg/ {2a \choose a} + \ldots \\
& + & {x^{a+1}y^{a-1}} + \ldots
\end{eqnarray*}
Thus:
\bq \beta_{H,-2} = {a^2-2a-1 \over 3a+1} {3a+1 \choose a}
\bigg/ {2a \choose a} \eq
and:
\bq R'_9(a,a,a)_{l,k} = {9a^2+6a+1 \over (3a+1)^2} {3a+1 \choose a}^2
\bigg/ {2a \choose a}^2 = {3a+1 \choose a}^2 \bigg/ {2a \choose a}^2 \eq
This is the desired answer.


\providecommand{\bysame}{\leavevmode\hbox to3em{\hrulefill}\thinspace}

\end{document}